\documentclass{article}

\usepackage{arxiv}

\usepackage[utf8]{inputenc} 
\usepackage[T1]{fontenc}    
\usepackage{hyperref}       
\usepackage{url}            
\usepackage{booktabs}       
\usepackage{amsfonts}       
\usepackage{nicefrac}       
\usepackage{microtype}      
\usepackage{lipsum}
\usepackage{graphicx}
\usepackage[ruled,vlined]{algorithm2e}
\usepackage{amsmath}
\usepackage{amssymb}
\usepackage{subfigure}
\usepackage{threeparttable}
\usepackage{cite} 
\usepackage{amssymb}

\newtheorem{theorem}{Theorem}

\newtheorem{assump}[theorem]{Assumption}

\title{Confounder Selection via Support Intersection}

\author{
  Shinyuu ~Lee \\
  \quad \\
  Peking University\\
  \texttt{ShinyuuLee@gmail.com} \\
   \And
 Yuru ~Zhu \\
 \quad\\
  Peking University\\
  \texttt{Zyr.stat@gmail.com} \\
}

\begin{document}
\maketitle

\begin{abstract}
Confounding matters in almost all observational studies that focus on causality.
In order to eliminate bias caused by connfounders, oftentimes a substantial number of features need to be collected in the analysis. In this case, large $p$ small $n$ problem can arise and dimensional reduction technique is required. However, the traditional variable selection methods which focus on prediction are problematic in this setting. Throughout this paper, we analyze this issue in detail and assume the sparsity of confounders which is different from the previous works. Under this assumption we propose several variable selection methods based on support intersection to pick out the confounders. Also we discussed the different approaches for estimation of causal effect and unconfoundedness test. To aid in our description, finally we provide numerical simulations to support our claims and compare to common heuristic methods, as well as applications on real dataset.
\end{abstract}


\section{Introduction}
\label{sec:headings}
Bias due to confounders  is one of the most important established systematic threats to consistent estimation of causal effects. Simpson's paradox is an elegant illustration of this type of bias that can arise in causal inference \cite{geng2019evaluation, pourhoseingholi2012control}. In order to identify causal effects in observational studies, researchers always rely on the unconfoundedness assumption, which requires that the individuals under different exposure level are in fact comparable once a sufficient set of pre-treatment covariates has been controlled for. 
%
In other literatures, unconfoundedness can be referred to using different terminology, such as exchangeability, weak ignorability or exogeneity \cite{vanderweele2019principles}.

To make the unconfoundedness assumption plausible, oftentimes a substantial number of features need to be collected in the analysis, but the estimated causal effect can be sensitive to the set of covariates included. More specifically, omitting confounding variables always lead to bias, on the other hand, including too many covariates may increase mean square error on the effect of interest and worsen the performance of standard estimation precedure, especially predictors only related to treatment, see \cite{greenland2008invited, ertefaie2018variable} for more detailed discussions. 
To obatin an unbiased, efficient estimation of causal effect, 
we could only pick out the confounders which affect both treatment assignment and outcome, since introducing them into the model is sufficient to remove the confounding bias without leading to inflation of variance. 
In this sense, selecting an appropriate set of confounders for which to control is critical for reliable causal inference, and sometimes useful to guide the trial design in the following phases to save time and financial resources. 

There have been numerous variable selection methods proposed in the last several decades,  such as the screening methods \cite{wang2009forward}, penalized methods LASSO and SCAD \cite{tibshirani1996regression, fan2001variable}, and Bayesian methods \cite{george2000calibration}. However, these variable selection methods focus on prediction, rather than estimating the effect of one exposure on response while treating the other predictors as confounders. 
Applying these general variable selection methods without specifically treating covariates as confounders could be problematic \cite{wilson2014confounder}. Confounder selection methods based on either just the treatment assignment model or just the outcome model may fail to account for non-ignorable confounders which barely predict the treatment or the outcome, respectively \cite{ertefaie2018variable}.
Belloni et al. (2014, 2017) show that attempting to control for high-dimensional confounders using a regularized regression adjustment obtained via, e.g., the lasso, can result in substantial biases\nocite{athey2018approximate, belloni2017program}. 

To address this problem, there has been considerable recent interest in performing confounder selection especially in high-dimensional settings.
Unfortunately, to the best of our knowledge, existing methods for confounder selection 
mainly requires the covariates which affect either treatment assignment or outcome are sparse, not confounders. A well-known representative approach is called double selection, 
removing the confounding bias by including the union of two sets of selected variables,
one selected by a lasso regressing the outcome on the covariates, and the other selected by a lasso logistic regression for the treatment assignment \cite{belloni2014inference}.
However, when the propensity or the outcome model is not sparse, we find that the performance of such double-selection methods is often poor.
Many other studies have such restrictives, for example, \cite{crainiceanu2008adjustment, wilson2014confounder, farrell2015robust, ma2019robust}. 
Recently, Athey, Imbens and Wager \cite{athey2018approximate} proposed a De-Biased technique and relax the sparse propensity conditions, but it still requires sparse regression model between covariates and outcome. 

In practice, when there are a large number of collected possible confounders, 
it is more reasonable to assume that many of them are related to treatment assignment or outcome, but the ones which affect both is sparse. While our goal should be identify these sparse confounders. 
Leveraging on these insights, we propose some feasible methods to perform confounder selection based on support intersection, only requring the sparsity of confounders.  Then we also discuss the approches to estimate the causal effect under our framework. 
We evaluate our method in the experiments using simulation datasets, as well as compare with existing methods. The results validate our theory and show the effectiveness of our method.

The rest of the paper is organized as follows.  Section 2 mainly includes some basic preliminaries, incuding our setting, assumptions as well as the notations.   In Section 3 we give the proposed confounder selection methods in high dimensions. Then  based on the selection  outcome, the estimation of causal effect and unconfoundedness test approaches is introduced in Section 4 and 5 respectively. At last, in Section 6 and 7 we provide our simulation and real data analysis results to support our claim.

\section{Notations and Assumptions}
\begin{figure}[h]
  \centering
  \includegraphics[width=0.5\textwidth]{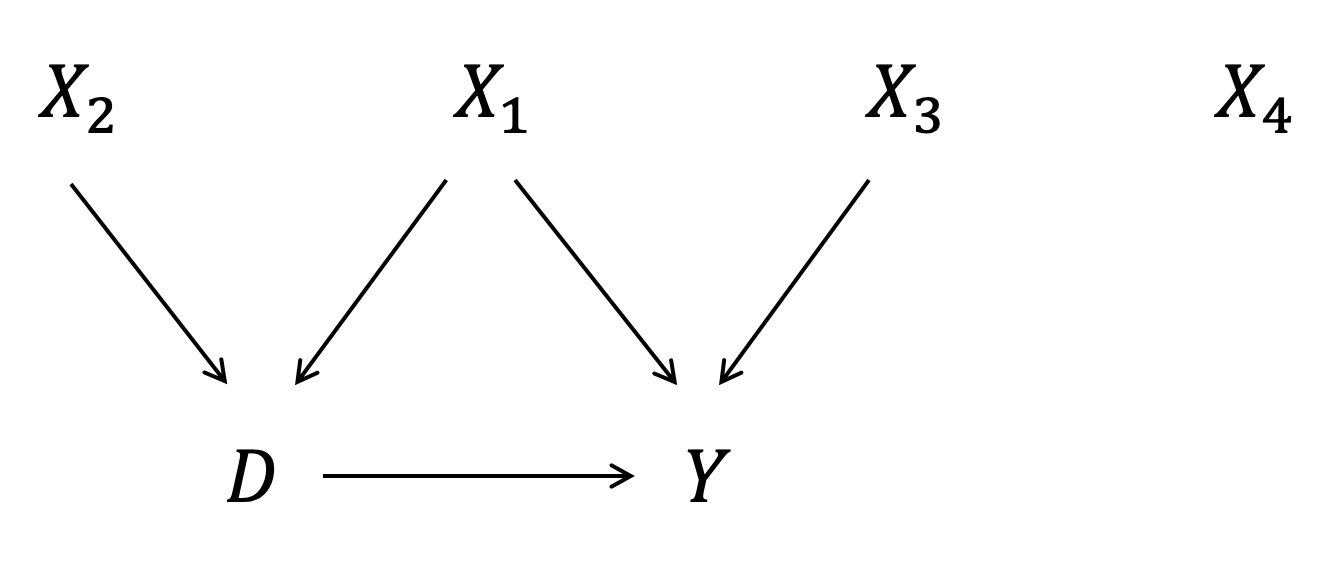}
  \caption{Covariate types}
  \label{fig:fig1}
\end{figure}
Following Neyman and Rubin \cite{rubin1974estimating}, we adopt the potential outcomes framework to define causal effects. Let $Y(d)$ denote the potential outcome that one would observe if the exposure  had, possibly contrary to fact, been set to level $d$, and let $D\in \mathbf{R}^{n}$ denote the treatment received. Let $X\in \mathbf{R}^{n\times p}$ denote the pre-treatment covariates which can be hign-dimentional.
Analogous to \cite{ertefaie2018variable},  we restrict attention here to the situation where each predictor can be classified into one of four following types
\begin{itemize}
    \item confounders (X1),which are related to both outcome and treatment;
    \item treatment predictors (X2),which are related to treatment and not to outcome; 
    \item outcome predictors (X3),which are related to outcome and not to treatment;
    \item redundant predictors (X4),which are not related to both outcome and treatment;
\end{itemize}{}
and the directed acyclic graph (DAG) in Figure \ref{fig:fig1} depicts the relationship intuitively.

The unconfoundedness assumption statistically states that the counterfactuals $Y(d)$ are independent of treatment $D$ conditional on some covariates $X$, which we denote by notation $Y(d) \perp D \mid X$. In this case, we say the covariates $X$ are suffice to control for confounding. Since we consider the confounder selection problem with a high-dimentional covariates, it is reasonable to assume unconfoundedness throughout this paper, in others words, all possible confounders have been included in collected covariates $X$. 
When the consistency assumption also natually holds, which states that $Y(d) =Y$  for those with $D=d$, then the causal effect is identifiable \cite{pourhoseingholi2012control, vanderweele2019principles}. 

\begin{assump}[sparsity]
We have a sequence of problems indexed by $n$, $p$, and $k$ such that the
number of confounders $X_1$ is no larger than $k$, and $k\ll \min(p,n)$.
\end{assump}

In this paper we consider the following linear model which is widely assumed in prior works, see \cite{athey2018approximate, farrell2015robust}. 
\begin{equation}\label{linear}
\begin{aligned}
 &   E(D\mid X;\alpha)  = g(X_1 \alpha_1 + X_2 \alpha_2) \\
 &   E(Y\mid D,X;\gamma,\beta)  = f(D\gamma + X_1 \beta_1 + X_3 \beta_2 )
\end{aligned}
\end{equation}
where $g$ and $f$ are some monotone functions. 
This linearity structure contains a  very wide range of models, regardless of whether $D$ and $Y$ are discrete or continuous. For example, when $Y$ is continuous and $D$ is binary, we let $f$ be identity function and $g$ be inverse logit function, so that the outcome and treatment models represent linear and logistic regression respectively. 
Without loss of generality, in the following section, we use the case where
$D$ and $Y$ are both continuous as representative with $f$ and $g$ being identity function, 
since this allows us to compare directly with traditional variable selection methods such as LASSO. In other cases, we just need to make the appropriate adjustments, like replacing according loss functions in the penalty methods.

Our goal is to distinguish the true confounders $X_1$ from other predictors, and ultimately give an unbiased, efficient estimation for $\gamma$, which represents the effect of $D$ on the response $Y$. For binary $D$, a more concerned statistic which could reflect the influence of $D$ on $Y$ is the average treatment effect (ATE) $E\{Y(1)-Y(0)\}$.

\section{Confounder Selection in High Dimensions}

In this section, we proposed two kinds of confounder selection algorithms, one is the screening method based on Pearson correlation, and another one is penalized method which turns our goal into an optimization problem. Both techniques are popular used in variable selection, and we   apply them to address our issue. It is assumed that the data has been centralized and standardized hereafter.

\subsection{Support Intersection SIS}

The sure independence screening (SIS) based on Pearson correlation learning is a method of reducing dimensionality from high to a moderate scale that is below the sample size. For a linear regression model $Y=X\beta+\epsilon$, after the design matrix $X$ is standardized columnwise so that their observed mean and sample variances equal to 0 and 1 respectively, the least square eatimator $\omega=\frac{1}{n}X^TY$ becomes a vector of the sample version of marginal Pearson correlations between individual predictor and the response variable rescaled by the standard deviation of the response. Fan and Lv \cite{fan2008sure} ranked the importance of features by ordering the components of $w$ decreasingly and selected the top covariates that are strongly correlated with the response. Specificly, for any give $\gamma\in(0,1)$, the submodel is 
\begin{equation*}
 \hat{\mathcal{M}}_r=\{1\leq i \leq p : |\omega_i|\text{ is among the first } [\gamma n]\text{ largest of all }\},
\end{equation*}

where $[\gamma n]$ denotes the integer part of $\gamma n$. This correlation screening procedure called SIS shrinks the full model down to the submodel with size $[\gamma n]$. SIS is invariant under scaling since only the order of componentwise magnitudes of $\omega$ is used. Note that $\frac{\lambda}{n}(X^TX+\lambda I_p)^{-1}X^TY\to\omega$ as $\lambda \to \infty$, which implies the componentwise regression estimator is a specific case of ridge regression with regularization parameter $\lambda =\infty$, namely, it reduces the variance of the estimator as much as possible. It is worth mentioning that the SIS merely ranks the importance of predictors, whereas the feature selection method based on multiple test judges the significance of each individual variable. Hence, further analysis of results of screening is needed. SIS has both low computational complexity and good property called the sure screening property. To be specific, assume that the true sparse model is $\mathcal{M}_*=\{1\leq j\leq p: \beta_j \neq 0\}$ with the size $s=|\mathcal{M}_*|<n$. Fan and Lv \cite{fan2008sure} proved that under four technical conditions, if $2\kappa + \tau < 1$ then there is some $\theta<1-2\kappa-\tau$ such that, when $\gamma \sim c n^{- \theta}$ with $c>0$, for some $C>0$, 
\begin{equation}
   P(\mathcal{M}_*\subset \hat{\mathcal{M}}_r)=1-O[\exp\{-C n^{1-2\kappa}/\log(n)\}]\to 1 
\end{equation}
as $n\to \infty$
and the true model size $s \leq [\gamma n]$ is implied by the assumptions. The property in (2) ensures all truly important variables survive after variable screening with probability tending to one. Any post-screening variable selection method is based on the screened submodels, thus the sure screening property is essential in practice. Although the above conditions are sometimes difficult to check, the numerical studies in \cite{fan2008sure} demonstrated that the SIS can efficiently reduce the ultrahigh dimension $p$ down to a relatively large scale $O(n^{1-\theta})$ for some $\theta>0$ and the submodel $\hat{\mathcal{M}}_r$ still  contains all important predictors with probability approaching one as $n$ tends to $\infty$ \cite{liu2015selective}. However, SIS may fail to select the important variables in some situations. Since SIS selects variables only according to the marginal correlations, it tends to select the unimportant predictor which is jointly uncorrelated but highly marginally correlated with the response, and an important predictor that is marginally uncorrelated but jointly correlated with the response may not be selected. Besides, the collinearity between predictors is another nonnegligible issure, some unimportant predictors that are highly correlated with the important predictors can have higher priority for being selected by SIS than other important predictors that are relatively weakly related to the response \cite{fan2008sure}. To enhance the performance of SIS, Fan and Lv \cite{fan2008sure} provided an iterative SIS procedure (ISIS). ISIS first uses an SIS-based model selection method, such as SIS-SCAD, to select a subset of $k_1$ variables $\hat{\mathcal{A}}_1=\{X_{i1}, \cdots,X_{i k_1}\}$  that were selected by SCAD on the basis of the joint information of remaining $[n/\log(n)]$  variables after correlation screening. Then iteratively treat the the the residual obtained from the regression of the responce on selected covariates in the previous step as the new response and apply the same SIS-based model selection method until the union of disjoint subsets obtained by each iteration has a desired size $d$ that is less than $n$.

Appling ISIS to confounder selection and based on the supprot intersection guidelines, we propose support intersection SIS, which can be implemented by the following procedures.

\begin{algorithm}
\caption{The Support Intersection SIS Algorithm}
\KwIn{The outcome $Y$, The treatment $D$, pre-treatment covariates $X$ and the maximum number of iterations $k$}
\KwOut{The set of selected confounders  $\hat{\mathcal{C}}=\hat{\mathcal{D}} \cap \hat{\mathcal{R}}$}
Initialize $\gamma$ as the corrsponding Ridge regression coefficient\\
Select the set of variables $\hat{\mathcal{D}}$  that are important for $D$ by ISIS.\\
\For{$i=1,\cdots,k$}
{
1. Select the set of variables $\hat{\mathcal{R}}$  that are important for $Y-D\hat{\gamma}$ by ISIS; \\
2. Treat variables in $\hat{\mathcal{C}}$ and $D$ as predictors and $Y$ as the response, use Ridge regression to update $\hat{\gamma}$; \\
3. \If{$\hat{\mathcal{R}}$ doesn't change} {\bfseries{break}}
\Else {\bfseries{Continue}}
}
\end{algorithm}

\subsection{Support Intersection LASSO}

Adding penalty term is a natural idea to perform simultaneous variable selection within treatment and outcome models. 
To solve our problem, we propose to use the $L_1$ norm of the Hadamard product (entrywise product) of $\alpha$ and $\beta$ as the penalty term. Essentially, this penalty term encourages the  sparsity of the support intersection  $\{ i:\alpha_i \neq 0 \} \cap \{i:\beta_i \neq 0 \}$.  In this way,  the optimization objective is 
\begin{equation}\label{obj0}
    \frac{1}{2} \parallel D- X\alpha \parallel_2^2 + \frac{1}{2} \parallel Y - D\gamma- X\beta \parallel_2^2 + \lambda  \parallel \alpha \circ \beta  \parallel_1
\end{equation}

For the $i$-th predictor, when $\alpha_i$ is punished to be zero or sufficiently small, then this penalty would not put restriction to $\beta_i$ any more, and vise versa.
In this sense, this penalty does not require the sparsity of either the support of the outcome model or treatment model, but only the sparsity of their intersection, therefore it meets our sparsity assumption. 

Now we consider some improvements to make our objective function more suitable for our setting. In (\ref{obj0}), the penalty term can be rewritten as $\sum_{i=1}^{p}\lambda |\alpha_i\beta_i|$, and the tuning parameter $\lambda$ always keeps same. This makes us lose modeling flexibility. Also, note that the penalty term can impose satisfying punishment on the $X_2$ and $X_3$, but far from enough on $X_4$ whose both corresponding coefficients are small. However, we expect the coefficients of $X_4$ near zero as much as possible. To this end, the corresponding $\lambda$ of $X_4$ should be larger than that of other predictors. All these reasons inspire us to design an adaptive tuning parameter,  putting relatively less punishment on the confounders and more punishment on the others.
An option is to let $\lambda_i \propto (\tilde{\alpha}_i \tilde{\beta}_i)^{-\frac{1}{2}}$, where $\tilde{\alpha}$ and $\tilde{\beta}$ are Ridge regression coefficients, and obviously it satisfies our requirement. 

Sometimes, one can add  $L_1$ norm of $\alpha$ and $\beta$ to avoid problems due to insufficient sample size, and accelerate the weak signal coefficients towards zero. 
In this case, we may only put a relative small tuning parameter and take this term as supplementary penalty. 
In the following section, unless otherwise specified, we adopt the following optimization objective to illustrate our idea, 
\begin{equation}\label{obj1}
    \frac{1}{2} \parallel D- X\alpha \parallel_2^2 + \frac{1}{2} \parallel Y - D\gamma- X\beta \parallel_2^2 + \sum_{i=1}^{p}\frac{\lambda}{\sqrt{ \tilde{\alpha}_i\tilde{\beta}_i}} |\alpha_i\beta_i| + \lambda_2 \parallel \alpha  \parallel_1 + \lambda_3 \parallel \beta  \parallel_1.
\end{equation}
The algorithm is given below where we apply the coordinate gradient descent to solve the convex optimization problem. As for how to choose hyperparameters, we can set a relatively small tuning parameter for the supplementary penalty, and use the following criterion analogous to GCV to choose the best $\lambda$,
    \begin{equation*}
        GCV^*(\lambda)= \frac{ \parallel D- X\hat{\alpha} \parallel_2^2 + \parallel Y - D\hat{\gamma}- X\hat{\beta} \parallel_2^2}{\left[1-(\parallel\alpha\parallel_0 +\parallel\beta\parallel_0)/2n\right]^2}.
    \end{equation*}{}

Related to our idea and method, Brandon et al.(2018) \cite{koch2018covariate} proposed a confounder selection procedure via adding penalty $\lambda \sum_{i}^{p} W_i\sqrt{\alpha_i^2+\beta_i^2}$, where weight $W_i$ can be adjusted adaptively to push relatively stronger punishment to covariates  unrelated to outcome. In this way, there are chances that some strong signals from treatment predictors would cover up the signal of confounders. That is one of the points why we choose to use the element-wise product instead of the sum of squares. 

\begin{algorithm}
\caption{The Support Intersection LASSO}
\KwIn{The outcome $Y$, treatment $D$, covariates $X$, the maximum number of iterations $k$, thershold $\epsilon$}
\KwOut{The estimated parameter $\hat{\alpha}$, $\hat{\beta}$ and $\hat{\gamma}$}
Initialize $\alpha$, $\beta$ and $\gamma$ as the corrsponding Ridge regression coefficients\\
\For{$i=1,\cdots,k$}
{
1. Update $\alpha$ by gradient descent; \\
2. Update $\beta$ by gradient descent; \\
3. Update $\gamma$ by minimizing the objective (\ref{obj1});\\
\If{change percentage of objective is smaller than $\epsilon$} {\bfseries{break}}
\Else {\bfseries{Continue}}
}
\end{algorithm}

\section{Estimation of Causal Effects}
\label{sec:others}

\subsection{Confounders Adjustment}

Once we pick out all the confounders among the predictors successfully, theoretically the confounding bias would disappear by introducing these variables into the model. If we adopt the selected confounders to perform adjustment, the key for consistency estimation is whether all the true confounders are identified. The output $\hat{\gamma}$ in SILASSO algorithm provides a pretty accurate estimator for $\gamma$, which can be verified in the simulations, and in the following we discuss some other techniques with well-established theoretical properties.

For simiplicity, we denote the selected confounders as $C$, whose output coefficient in either outcome model or treatment model is away from zero.
When both $Y$ and $D$ are continuous and the linearity holds, a simple way is to rebuild the conditional outcome mean model $E(Y|D,C)$ and take the coefficient of $D$ as the estimator for $\gamma$. This is a relatively low dimensional regression problem. Moreover, we could also consider introduce the selected outcome predictors into the model, and if most of them are indeed the true outcome predictors, it would promote the model explanatory ability, and in this case, this approach can be a good choice for us.

 It is worth noting that when $D$ is binary, the average causal effect (ATE) can be doubly robust estimated \cite{bang2005doubly} by 
\begin{equation*}
\small
\mathbf{P_n}\left[  
\frac{DY}{E(D| X;\hat{\alpha})}
-\frac{(1-D)Y}{1-E(D|X;\hat{\alpha})}
- \frac{D-E(D| X;\hat{\alpha})}{E(D| X;\hat{\alpha})}E(Y| 1,X;\hat{\gamma},\hat{\beta})
- \frac{D-E(D| X;\hat{\alpha})}{1-E(D| X;\hat{\alpha})}E(Y| 0,X;\hat{\gamma},\hat{\beta})
\right],
\end{equation*}
and we could substitute the covariates $X$ by the corresponding selected predictors.

\subsection{Instrument Variables Adjustment}

In our setting, the treatment predictors $X_2$ satisfy the following three conditions: (i) they have causal effect on $D$, (ii) only affect the response $Y$ through the treatment $D$, and (iii) they do not share common causes with the outcome $Y$. Therefore the treatment predictors are actually instrument variables \cite{hernan2006instruments}.   

This inspires us to utilize these  instrument variables to estimate the causal effects, which can be less sensitive to unconfoundedness assumption. In this sense, even if we omit some confounders, as long as we can pick out several instrument variables, then we can still identify and estimate the causal effects consistently. Both the previously mentioned SISIS and SILASSO methods could help us pick out potential instrumental variables. 
In SILASSO, we select the ones whose output coefficient in outcome model is zero while in treatment model is not zero. In what follows, we denote the selected instrument variables as $Z$.

When there exist multiple instrument variables,  there have been many mature methods to address this problem. A common class of such scenario in causal inference is Mendelian randomization, which use large numbers of genetic variants as instrumental variables in observational data \cite{burgess2017review}. With multiple instrument variables, a popular way to estimate the causal effect of the exposure on the outcome is the two-stage least squares (TSLS) method. However, when some instrument variables are invalid, the estimator derived by this method can be biased. Jack Bowden et al.(2015) \nocite{bowden2015mendelian} proposed a robust estimator through Egger regression, allowing us estimate the causal effect consistently under a weaker condition called InSIDE (Instrument Strength Independent of Direct Effect).  In general, these methods basically require the independence between instrument variables, as well as strong validness. 

In our setting, the instrument variables may be related to each other, and moreover, we could not be completely convinced that the selected ones must be instrument variables. 
Therefore we consider another robust method which requires less restrictions, propoed by Frank Windmeijer et al. (2019)\nocite{windmeijer2019use}. When the true instrument variables accounts for more than half of the total selected instrument variables, under some regular conditions, they have shown that the median of the coefficients ratio $\hat{\pi}$ is a consistent estimator for $\gamma$, where $\hat{\pi}_j= {\hat{\Gamma}_j}/{\hat{\eta}_j}$, $\hat{\Gamma} = (Z^{\prime}Z)^{-1}Z^{\prime}Y$ and $\hat{\eta} = (Z^{\prime}Z)^{-1}Z^{\prime}D$.
However this median converges in distribution to that of an
order statistic, not a normal distribution. To remedy this defect, they apply the adaptive Lasso method  and give an asymptotically normal and consistent estimator for $\gamma$ as 
\begin{equation*}
    \frac{\hat{D}^{\prime}(Y- Z\hat{\alpha}_{ad})}{\hat{D}^{\prime}\hat{D}} \quad  \hat{\alpha}_{ad} = \arg \min  \frac{1}{2}\parallel Y- \tilde{Z}\alpha \parallel_2^2 +\lambda\sum_{l}\frac{|\alpha_l|}{|\hat{\alpha}_{m,l}|^v}
\end{equation*}{}
where $\hat{D}=Z\hat{\eta}$, $\hat{\alpha}_m = \hat{\Gamma} - \mbox{median}(\hat{\pi}) \hat{\eta}$, $\tilde{Z}=(I-\hat{D}(\hat{D}^{\prime}\hat{D})^{-1}\hat{D}^{\prime})Z$.

In the following simulation, we would assess the performance of our algorithm to see whether the true ones can account for more than half
among the instrument variables selected, which is the key for consistency estimation.

\section{Testing Unconfoundedness}

\paragraph{Using Instrumental Variables}
As we have discussed in the previous section, the treatment predictors can be seemed as instrument variables. With the help of instrument variable, one can test the unconfoundedness assumption. The existing and widely used method for testing whether there exist unobserved confounder is the Durbin-Wu-Hausman endogeneity test (DWH test), independently proposed by Dur bin, Wu and Hausman. This test only needs one instrument variable, but it requires it is a valid one. Note DWH test assumes that the treatment effect is homogeneous, therefore if this test rejects, one cannot be sure whether it is because of unmeasured confounding or treatment effect heterogeneity. For binary $D$ and binary instrument variable,  Donald, Hsu and Lieli (2014) proposed such a Durbin-Wu-Hausman type statistic \cite{donald2014testing}.

\paragraph{Using Outcome Predictors}
If we could pick out some outcome predictors, then these variables could also help us perform test on unconfoundedness. In SILASSO, we can select the ones whose output coefficient in treatment model is zero while in outcome model is not zero. Such a procedure is proposed by Zongwu Cai et al.\ (2019) to test unconfoundedness and it relies on the existence of an auxiliary variable which is correlated to potential outcomes but is independent of the treatment status given on potential outcomes and observable covariates. Obviously, any outcome predictor is satisfied the requirement. See \cite{cai2019testing} for more details.

\paragraph{Sensitivity Analysis}
When we  are not sure whether there are unobserved confounders, sensitivity analysis provides us a way to assess the influence of potential confounding bias on the estimated causal effect. For example, the derived conclusion shows there exist a positive effect from $D$ on $Y$, but this may be an illusion caused by confounding bias, while sensitivity analysis can tell us how strong the unobserved confounder should be to have such ability.
If we have confidence that there is no such powerful confounder, then our original conclusion is valid. The detailed description is available in \cite{robins2000sensitivity}.

\section{Simulation}
\subsection{Design}
We investigate the performance of SISIS and SILasso in different scenarios and compare it with the results of other approaches to selecting confounders. To illustrate the numerous simulation specificly, we introduce some notations. The dimension of $X_1, X_2, X_3, X_4$ is denoted as $p_1, p_2, p_3, p_4$ respectively. For each scenario, we generate $n$ observations independently from $p_i$ dimensional normal distribution $\mathcal{N}_{p_i}(0,\Sigma_i)$ to form the data set of $X_i$, where the element at the $j$-th row and $k$-column of the covariance matrix $\Sigma_i$ is $\sigma^2\rho^
{|j-k|}, i=1,\cdots,4$. Besides, the true coefficients of predictors, i.e. $\alpha_1,\beta_1$, $\alpha_2,\beta_2$, are generated randomly from the uniform distribution with the range $(-1,-0.2)\cup(0.2,1)$. Additionally, each component of the error terms $\epsilon_i$ is generated from $N(0,\sigma^2),i=1,2$. The value of $\gamma$ can be set to 1 without the loss of generality. Then, the data generating mechanism of the treatment and response can be summarized as $D=X_1 \alpha_1 + X_2 \alpha_2+\epsilon_1$ and $ Y=D + X_1 \beta_1 + X_3 \beta_2+\epsilon_2$.

For the $i$-th predictor, let $L_{i}^{m}=|\hat{\alpha}_i^m \hat{\beta}_i^m|$ for the model $m$. Note that 
supposed we know the real $\alpha_i$ and $\beta_i$ (which is impossible in practice), only the $|\alpha_i\beta_i|$ of confounders  should differ from zero. Therefore based on $L_{i}^{m}$, we can derive an importance rank for all the selected variables, and a predictor with larger $L_{i}^{m}$ is more likely to be regarded as the confounder. Leveraging on these insights, we use the largest ordinal number of the true confounders in the selected sequence as the our evaluation tool  to assess the performance of each model, denoted as \textit{Cover Number  (CN)}.

We consider three scenarios: (1) $n=1200, p_1=10, p_2=p_3=300, p_4=390, \sigma^2=1, \rho=0.2$; (2) $n=1200, p_1=10, p_2=p_3=300, p_4=390, \sigma^2=1, \rho=0.4$; (3) $n=1200, p_1=10, p_2=p_3=500, p_4=490, \sigma^2=1, \rho=0.2$ and compare the performance of four methods: ordinary least square (OLS), Double Selection via LASSO (LASSO hereafter), Support Intersection SIS (SISIS) and Support Intersection LASSO (SILASSO).
\subsection{Results}
Table \ref{tab:table} shows the number of selected confounders (SN),  the number of true confounders selected (TN) and the corresponding Cover Number (CN) by four methods respectively. We can see that under the situation where $p<n$, compared to OLS, other three methods all enhance the ability to discriminate confounders more or less. At the same time, the SISIS and SILASSO select much less pontential confounders including all ten true confounders, thus they perform better than the other two ordinary methods. When the correlation between predictors are larger, it can be found that the new proposed methods are still quite stable. As for the high dimensional large $p$ small $n$ setting, the SISIS only selects three true confounders, which is dissatisfactory. While SILASSO always provides a good and stable option. 

 \begin{table}[h!]
 \caption{Comparison of Four Methods}
  \centering
  \renewcommand\tabcolsep{9.0pt} 
  \begin{threeparttable}
  \begin{tabular}{llllllllll}
    \toprule
    Methods    &  \multicolumn{3}{c}{ Scenario 1}    & \multicolumn{3}{c}{ Scenario 2}  & \multicolumn{3}{c}{ Scenario 3}   \\
    \midrule
    & SN & TN & CN & SN & TN & CN &SN & TN & CN\\
        \cmidrule(r){2-10}
    OLS & --- & --- & 122  & --- &--- & 57  &  ---&---& 800 \\
    LASSO   & 495 & 10  & 21 &  546 & 10 &38 &  668 & 7 & --- \\
    SISIS    & 17 & 10 & 14 & 45 &10 & 45    &11 & 3 & --- \\
    SILASSO  & 23 & 10&10  & 29 &10 & 10 & 91 & 10 & 10\\
    \bottomrule
  \end{tabular}
  \label{tab:table}
  
   \begin{tablenotes}
        \footnotesize
        \item[1] SN, TN and CN are short for the number of selected variables as confounders, the number of true confounders among selected variables and Cover number respectively.  
        \item[2] If the method successfully selects all confounders, TN equals 10 and the closer CN is to TN, the stronger ability to discriminate confounders. 
        \item[3] TN $<10$ indicates that the method is failed and therefore CN is meaningless.
        \item[4] Sinc OLS can not perform confounder selection, its SN and TN do not exist.
       \end{tablenotes}
    \end{threeparttable}
\end{table}

\begin{figure}[h!]
  \centering
  \includegraphics[width=0.45\textwidth]{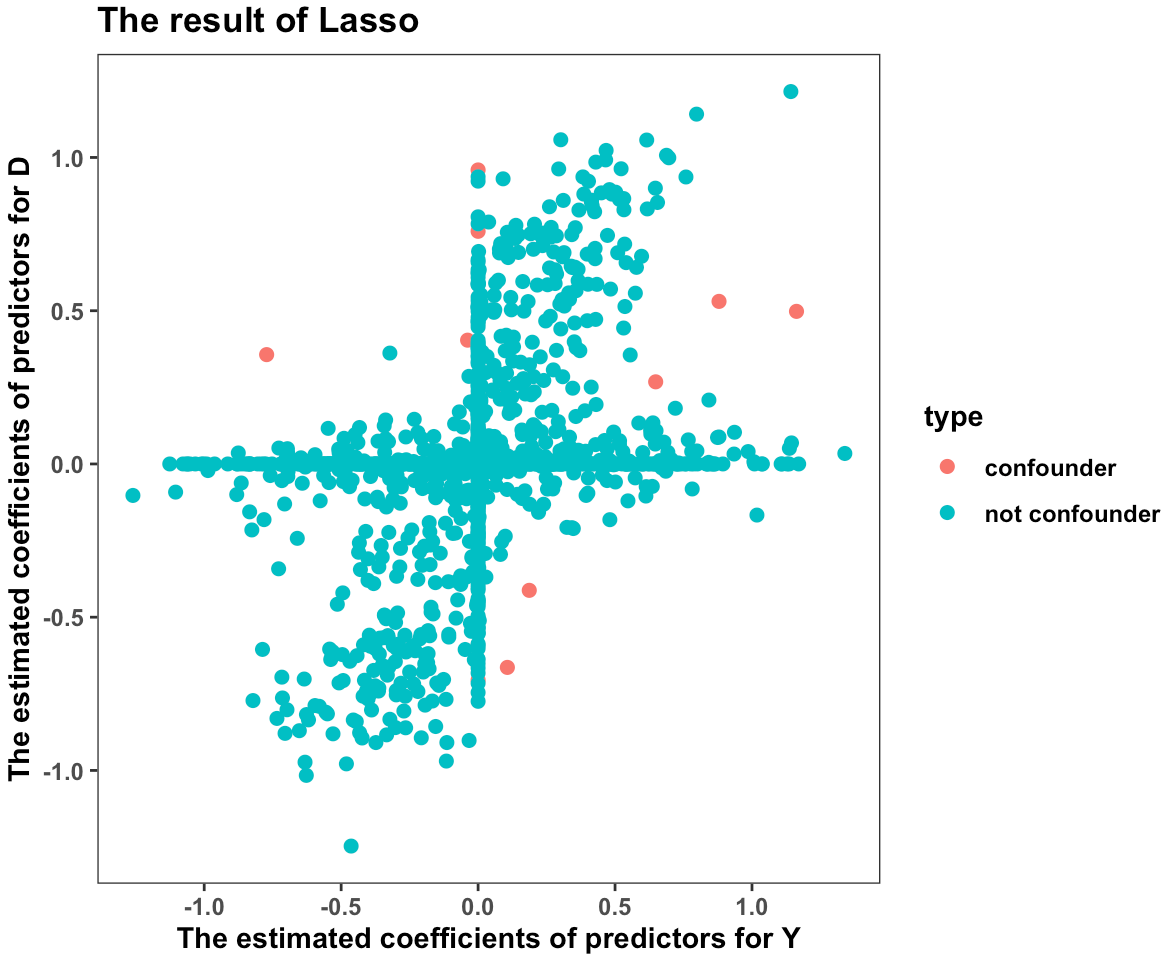}
  \includegraphics[width=0.45\textwidth]{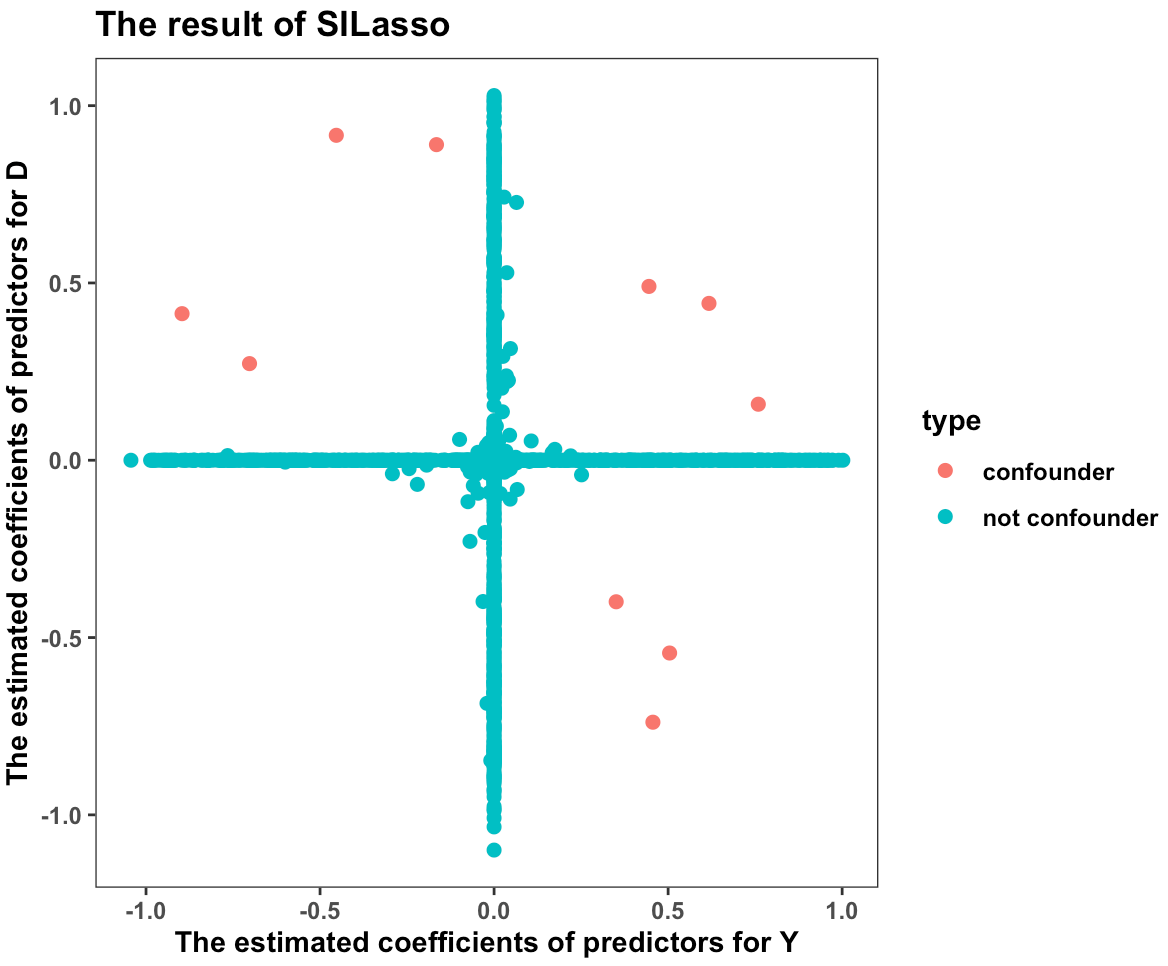}
  \caption{\textbf{The results of LASSO and SILASSO (p=1500, n=1200)} When a point is punished to lie on the X or Y axis, it indicates that the corresponding variable only affects $Y$ or $D$. On the other hand, a point which is far away from the coordinate axes can be regarded as the confounders affecting both $D$ and $Y$. The red points are the true confounders, and it can be seen that SILASSO successfully discriminates them while LASSO fails.   }
  \label{fig:fig2}
\end{figure}

The estimated coefficients of $X$ obtained from the SILASSO method in scenario 3 are ploted in figure \ref{fig:fig2}. The points corresponding to true confounders are red, which are far away from the coordinate system. It can be clearly seen on this graph that, except for a small number of points, most of the points are forced to lie  on the X or Y axis. This helps us to identify confounders quickly, which supports the superiority of our proposed method.

\section{Application on Real Data}

\subsection{Data Description} 
We examine the performance of our proposed methods on the determinants of economic growth dataset firstly used by Sala-i-Martin et al.(2004) \cite{sala2004determinants},
which contains observations for 139 countries and 68 variables. Following Ashkan et al. (2015) \cite{ertefaie2018variable}, we consider a subset of the data which includes 88 countries and 35 variables for illustration purposes, and are interested in selecting the non-ignorable covariates which confound the effect of \textit{public education spending share of GDP} (exposure) as a measure of investment in education on \textit{the average growth rate of GDP per capita in 1960-1996} (outcome).
This data is described in detail by Gernot Doppelhofer et al. (2009) \cite{doppelhofer2009jointness}.

\subsection{Analysis Results}

We apply our methods to the determinants of economic growth dataset which is described previously, selecting the confounders which affect both the exposure (\textit{public education spending share of GDP})  and outcome (\textit{the average growth rate of GDP per capita in 1960-1996}).
The SISIS selects six confounder, while the SIlASSO selects five confounders, and the results are shown in Table \ref{tab:table2}. It can be seen that the results of two methods are similar. The selected confounders are reasonable, for example, the fraction of time spent in war can both affect public education spending and the economy growth.

 \begin{table}[h!]
 \caption{The Results on Real Data Set }
  \centering
  \renewcommand\tabcolsep{9.0pt} 

  \begin{tabular}{ccc}
    \toprule
   Name of Selected Confounder      &  SISIS    &  SILASSO    \\
    \midrule 
    Fraction of population Catholics  & \checkmark &  \checkmark \\
    Average investment price level on purchasing power parity basis   & \checkmark & \checkmark  \\
    Enrolment rate in primary education     & \checkmark & \checkmark  \\
   Fraction of time spent in war & \checkmark & \checkmark\\ 
       Coastal population per coastal area & & \checkmark\\
    Fraction of population Protestant  & \checkmark &  \\
    Number of years economy has been open  &\checkmark & \\
    \bottomrule
  \end{tabular}
  \label{tab:table2}
  
\end{table}

We introduce the selected confounders into the model and estimate the corresponding coefficient of exposure as the  causal effect estimator. The p-value shows that the effect is not significant at a confidence level of 0.05, which is consistent with the result we obtained under the full model including all variables. The conclusion looks counterintuitive, as generally it is believed that investment in education can help promote the speed of economic growth. However, the actual situation may be more complicated, and many literatures have confirmed this statement. We quote the following words from the document \textit{Economic Returns to Investment in Education} published by the world bank \cite{galal2007road},  

\small

\quad \textit{But does a higher level of investment in education affect the growth path? The answer to the latter question is predominantly “no.” Barro and Lee (1994) show that the increase in the number of those who attended secondary school between 1965 and 1985 had a positive effect on growth, but estimates by others do not confirm this result. Using an aggregated production function, Benhabib and Spiegel (1994) and Pritchett (1996) also measure the impact of human capital investment on the rate of economic growth. They use various measurements of human capital, including the number of years of education, literacy rates, and secondary enrolment rates. Whatever the education variable chosen, the associated coefficients appear either as insignificant or as having a negative sign.} 

\normalsize
A possible reason given by this document is that the countries which had a higher level of education in 1960 had a greater opportunity, 40 years later, to reach a higher level of development.  In other words, the impact of investment in education on economic growth has a long-term lag effect. Of course, in view of that this article focuses on confounder selection, the relationship between investment in education and economy growth is  another story left to economists. 

\bibliographystyle{unsrt}  
\bibliography{references}  

%

\end{document}